\documentclass[12pt]{article}  
\usepackage{amsmath}
\usepackage{amssymb}

\newtheorem{theorem}{Theorem}

\newcommand{\CC}{{\Bbb C}}

\newcommand{\PP}{{\Bbb P}}

\newcommand{\ZZ}{{\Bbb Z}}

\newcommand{\eps}{\varepsilon}

\title{Indices of vector fields or 1-forms and characteristic numbers}
\author{W.Ebeling and S.M.Gusein-Zade
\thanks{Partially supported by the DFG-programme ''Global methods in
complex geometry'', grants RFBR--01--01--00739, INTAS--00-259,  NWO--RFBR--047.008.005.
Keywords: characteristic numbers, 1-forms, vector fields, singular points, indices.
AMS Math. Subject Classification: 32S50, 57R20, 57R25.
}
}
\date{}

\begin{document}

\maketitle

\begin{abstract} We define an index of a collection of 1-forms on a complex isolated
complete intersection singularity corresponding to a Chern number and, in the case when the
1-forms are complex analytic, express it as the dimension of a certain algebra.
\end{abstract}

\section*{Introduction}
An isolated singular point (zero) $p$ of a vector field $X$ on a
(real or complex analytic) manifold $M$ has an index: the degree of the map
$X/\|X\|$ from a small sphere around the point $p$ to the unit sphere. If
the manifold $M$ is closed (compact without boundary) and the vector field
$X$ has only isolated singular points, the sum of their indices is equal to
the Euler characteristic $\chi(M)$ of the manifold $M$. If $M$ is a complex
analytic manifold of dimension $n$, then its Euler characteristic $\chi(M)$
is the characteristic number $\langle c_n(TM), [M]\rangle$, where $TM$ is
the tangent bundle of the manifold $M$ and $c_n$ is the corresponding Chern
class.

If the vector field $X$ is complex analytic, the index of its (isolated)
singular point $p$ can be computed as the dimension of a certain algebra.
Namely, if, in local coordinates centred at the point $p$, $X= \sum_{i=1}^n
X_i (\partial/\partial x_i)$, then $\mbox{ind}_p X = \dim {\cal
O}_{\CC^n,0}/(X_1, \ldots , X_n)$, where ${\cal O}_{\CC^n,0}$ is the ring of
germs of holomorphic functions of $n$ variables, $(X_1, \ldots, X_n)$ is the
ideal generated by the components of the field $X$. The
Eisenbud-Levine-Khimshiashvili formula \cite{EL, Kh} expresses the index of an
algebraically isolated singular point of a (smooth) vector field on a real
smooth manifold as the signature of a certain quadratic form defined on the real
analogue of the algebra mentioned above. Similar statements hold for 1-forms on
manifolds with the only difference that the sum of indices of singular points of
a 1-form on a complex analytic manifold $M$ of dimension $n$ is equal to the
characteristic number
$\langle c_n(T^\ast M), [M]\rangle = (-1)^n \chi(M)$.

In the framework of the described properties there is almost no difference
between vector fields and 1-forms (except the fact that non-zero analytic
vector fields and 1-forms exist, generally speaking, on different complex
manifolds). The difference becomes crucial if one considers complex-analytic
vector fields or 1-forms on singular varieties \cite{EG1, EG2}. The main
difference is related with the fact that a complex analytic 1-form can be
restricted to a submanifold or to a subvariety (staying complex analytic, at
least, at the non-singular points of the subvariety), while there is no
invariant notion of a restriction of a vector field to a submanifold.

Here we start to discuss generalizations of these notions and
statements for other characteristic numbers (different from $\langle
c_n, [M]\rangle$). For smooth varieties, it is not very essential whether we
discuss vector fields, or 1-forms, or simply sections of a vector bundle. For
singular varieties (namely for isolated complete intersection singularities: ICIS) the
properties of 1-forms and vector fields with respect to the discussed problems are quite
different. Here we restrict ourselves almost exclusively to the complex analytic case.

\section{An index for collections of sections of a vector bundle}

We consider a general setting. Let $\pi : E \to M$ be a complex analytic
vector bundle of rank $m$ over a complex analytic manifold $M$ of dimension
$n$. It is known that the Poincar\'e dual to the characteristic class
$c_k(E)$ ($k=1, \ldots , m$) is represented by the $2(n-k)$-dimensional cycle which consists
of points of the manifold $M$ where $m-k+1$ sections of the vector bundle
$E$ are linearly dependent (cf., e.g., \cite[p.~413]{GH}).

For natural numbers $p$ and $q$ with $p \geq q$, let ${\cal M}(p,q)$ be
the space of $p \times q$ matrices with complex entries and let $D_{p,q}$
be the subspace of ${\cal M}(p,q)$ consisting of matrices of rank less than
$q$.
The subset $D_{p,q}$ is an irreducible subvariety of
${\cal M}(p,q)$ of codimension $p-q+1$ (see below). The complement $W_{p,q} =
{\cal M}(p,q) \setminus D_{p,q}$ is the Stiefel manifold of
$q$-frames (collections of $q$ linearly independent vectors) in $\CC^p$.
It is known that $W_{p,q}$ is $(2p-2q)$-connected and $H_{2p-2q+1}(W_{p,q})
\cong \ZZ$ (see, e.g., \cite{H}). The latter fact also proves that $D_{p,q}$ is irreducible.
Since $W_{p,q}$ is the complement of an irreducible complex analytic subvariety of
codimension $p-q+1$ in ${\cal M}(p,q)$, there is a natural choice of a
generator of the homology group $H_{2p-2q+1}(W_{p,q}) \cong \ZZ$. Namely, the
(''positive'') generator is the boundary of a small ball in a smooth complex
analytic slice transversal to $W_{p,q}$ at a non-singular point. 

Let ${\bf
k}=(k_1, \ldots , k_s)$ be a sequence of positive integers with
$\sum_{i=1}^s k_i = k$.
Consider the space ${\cal M}_{m, {\bf k}}=
\prod_{i=1}^s {\cal M}(m,m-k_i+1)$ and the subvariety $D_{m, {\bf k}}=
\prod_{i=1}^s D_{m,m-k_i+1}$ in it. The variety $D_{m, {\bf k}}$ consists of
sets $\{A_i\}$ of $m \times (m-k_i+1)$ matrices such that $\mbox{rk}\, A_i <
m-k_i+1$  for each
$i=1, \ldots , s$. Since $D_{m, {\bf k}}$ is irreducible of codimension
$k$, its complement $W_{m, {\bf k}}= {\cal M}_{m, {\bf k}} \setminus D_{m, {\bf k}}$ 
is $(2k-2)$-connected,
$H_{2k-1}(W_{m, {\bf k}}) \cong \ZZ$, and there is a natural choice of a
generator of the latter group. This choice defines a degree (an integer) of
a map from an oriented manifold of dimension $2k-1$ to the manifold $W_{m, {\bf k}}$.

Let $\{\omega^{(i)}_j\}$ ($i=1, \ldots , s$; $j=1, \ldots , m-k_i+1$;
$\sum_{i=1}^s k_i = n$) be a
collection of sections of the vector bundle $\pi:E \to M$ such that there
are only isolated points where, for each $i=1, \ldots, s$, the sections
$\omega^{(i)}_j$, $j=1, \ldots , m-k_i+1$, are linearly dependent. Let us
choose a trivialization of the vector bundle $\pi: E \to M$ in a
neighbourhood of 
a point $p$, let $(\omega^{(i)}_1(x), \ldots, \omega^{(i)}_{m-k_i+1}(x))$ be the $m
\times (m-k_i+1)$-matrix the columns of which consist of the components of the sections
$\omega^{(i)}_j(x)$, $j=1, \ldots , m-k_i+1$, $x \in M$, with respect to this
trivialization. 
Let
$\Psi_p$ be the mapping from a neighbourhood of the point $p$ to ${\cal M}_{m, {\bf k}}$
which sends a point $x$
to the collection of matrices $\{ (\omega^{(i)}_1(x), \ldots,
\omega^{(i)}_{m-k_i+1}(x))\}$, $i=1, \ldots, s$.
Its restriction
$\psi_p$ to a small sphere
$S^{2n-1}_{p,
\eps}$ around the point $p$  maps this sphere to the subset $W_{m, {\bf k}}$. 

\addvspace{3mm}

\noindent {\bf Definition} The degree of the mapping $\psi_p :
S^{2n-1}_{p, \eps} \to W_{m, {\bf k}}$ is called the {\em index} of the
collection of sections $\{\omega^{(i)}_j\}$ at the point $p$ and is denoted by
$\mbox{ind}_p\{\omega^{(i)}_j\}$.

\addvspace{3mm}

In other words, the index $\mbox{ind}_p\{\omega^{(i)}_j\}$ is equal to the intersection
number of the germ of the image of the (complex analytic) map $\Psi_p$ with the variety
$D_{m, {\bf k}}$. 

The point $p$ is non-singular for the collection $\{\omega^{(i)}_j\}$ if at least
for some $i$ the values
$\omega^{(i)}_j(p)$ of the sections $\omega^{(i)}_j$ at the point $p$ are
linearly independent. The index $\mbox{ind}_p \{\omega^{(i)}_j \}$ of a non-singular point
is equal to zero. The description of the Poincar\'e dual cycles to the Chern classes
$c_{k_i}(E)$ of the vector bundle $E$ yields the following statement.

\begin{theorem} \label{Th1}
Let $\sum_{i=1}^s k_i  = n$ and
suppose that the collection $\{\omega^{(i)}_j\}$ ($i=1, \ldots , s$; $j=1, \ldots,
m-k_i+1$) of sections of the vector bundle $\pi: E
\to M$ over a closed complex manifold $M$  has
only isolated singular points. Then the sum of the indices of these points is equal to the
characteristic number
$\langle
\prod_{i=1}^s c_{k_i}(E), [M]\rangle$ of the vector bundle $E$.
\end{theorem}

\section{An algebraic formula for the index}
Suppose that, in local coordinates centred at the point $p$, the components of the
section $\omega^{(i)}_j$ are complex analytic functions on $M$ (e.g., this
happens if the fibre bundle $\pi : E \to M$ is complex analytic and the
sections $\omega^{(i)}_j$ are complex analytic in a neighbourhood of the point
$p=0$). Let $I_{\{\omega^{(i)}_j\}}$ be the ideal in the ring ${\cal
O}_{\CC^n,0}$ of germs of analytic functions of $n$ variables generated by
the $(m-k_i+1) \times (m-k_i+1)$-minors of the matrices
$(\omega^{(i)}_1, \ldots, \omega^{(i)}_{m-k_i+1})$ for all
$i=1, \ldots, s$.

\begin{theorem} \label{Th2}
$${\rm ind}_p\{\omega^{(i)}_j\} = \dim_\CC {\cal
O}_{\CC^n,0}/I_{\{\omega^{(i)}_j\}}.$$
\end{theorem}

\noindent {\em Proof.} The proof follows the same lines as the proof of
\cite[Theorem~1]{EG2}. Let $\widetilde{\omega}^{(i)}_j$ be suitable analytic
perturbations of the sections $\omega^{(i)}_j$ and let $\widetilde{D}$ be the
corresponding degeneracy locus. The essential point in the proof is the fact
that $\widetilde{D}$ is Cohen-Macaulay which follows from
\cite[Corollary~(2.7)]{BR64}. $\Box$

\section{Indices on ICIS}
Let $f=(f_1, \ldots , f_\ell) : (\CC^n,0) \to (\CC^\ell,0)$ be an analytic map which
defines an $(n-\ell)$-dimensional isolated complete intersection singularity (ICIS)
$V=f^{-1}(0)
\subset (\CC^n,0)$ ($f_r:(\CC^n,0) \to (\CC,0)$). Let $\{ \omega^{(i)}_j\}$ be a
collection of 1-forms on a neighbourhood of the origin in
$(\CC^n,0)$ with $i=1, \ldots, s$, $j=1, \ldots , n-\ell-k_i+1$, $\sum k_i = n-\ell$,
the restriction of which to $V$ has no singular points on $V$ outside of the origin in a
neighbourhood of it. (Here it is not necessary to demand that the 1-forms $\omega^{(i)}_j$
are complex analytic. It is sufficient to suppose that $\omega^{(i)}_j$ are continuous
complex linear functions on the tangent bundle $T\CC^n$.) Let
$U$ be a neighbourhood of the origin in
$\CC^n$ where all the functions
$f_r$ ($r=1, \ldots , \ell$) and the 1-forms $\omega_j^{(i)}$ are defined and such
that the restriction of the collection $\{ \omega_j^{(i)} \}$ of 1-forms to $(V \cap
U) \setminus \{ 0\}$ has no singular points. Let $S_\delta \subset U$ be a
sufficiently small sphere around the origin which intersects $V$ transversally
and denote by $K=V \cap S_\delta$ the link of the ICIS $(V,0)$. Let $\widehat{\bf k}
= (k_1+\ell, \ldots , k_s +\ell)$ and let $\Psi_V$ be the mapping from $V \cap U$ to
${\cal M}_{n,\widehat{\bf k}}$ which sends a point $x \in V \cap U$ to the collection of $n
\times (n-k_i+1)$-matrices 
$$\{ (df_1(x), \ldots , df_\ell(x), \omega_1^{(i)}(x), \ldots ,
\omega_{n- \ell -k_i+1}^{(i)}(x)) \}, \quad i=1, \ldots, s.$$
Its restriction $\psi_V$ to the link $K$ maps $K$ to the subset
$W_{n,\widehat{\bf k}}$. 

\addvspace{3mm}

\noindent {\bf Definition} The {\em index} ${\rm ind}_{V,0}\{\omega^{(i)}_j\}$ of the
collection of 1-forms $\{ \omega^{(i)}_j \}$ on the ICIS $V$ is the degree of the mapping
$\psi_V : K \to W_{n, \widehat{\bf k}}$. 

\addvspace{3mm}

As above, the index ${\rm ind}_{V,0}\{\omega^{(i)}_j\}$ is equal to the intersection
number of the germ of the image of the complex analytic mapping $\Psi_V$ with the variety
$D_{n, \widehat{\bf k}}$.

For $s=1$, $k_1=n-\ell$, this definition coincides with the one of the index of a 1-form
from \cite{EG1}.

Let $V \subset \CC\PP^n$ be a complete intersection with isolated singular
points. Let $L$ be a complex line bundle on $V$ and let
$\{ \omega^{(i)}_j \}$ be a collection of continuous 1-forms on $V$ with values
in $L$. Here this means that the forms $\omega_j^{(i)}$ are continuous
sections of the fibre bundle $T^\ast V \otimes L$ outside of the singular points
of $V$. 
Since, in a neighbourhood of each point $p$, the vector bundle $L$ is trivial,
one can define the index $\mbox{ind}_p \{\omega^{(i)}_j\}$ of the collection of
1-forms $\{ \omega_j^{(i)}\}$ at the point $p$ just in the same way as in the
local setting above. Moreover, one can consider $L$ as a line bundle on a smoothing
$\widetilde{V}$ of the complete intersection $V$ as well (e.g., using the pull back along a
projection of $\widetilde{V}$ to $V$).

\begin{theorem} \label{Th3}
One has
$$ \sum_{p \in V} {\rm ind}_p \{\omega^{(i)}_j\} = \langle \prod_{i=1}^s
c_{k_i}(T^\ast\widetilde{V} \otimes L), [ \widetilde{V} ] \rangle,$$
where $\widetilde{V}$ is a smoothing of the complete intersection $V$.
\end{theorem}

Now, as above, let $(V,0) \subset (\CC^n,0)$ be the ICIS defined by the equations $f_1=
\cdots = f_\ell=0$ and let $\{ \omega_j^{(i)}\}$ ($i=1, \ldots, s$; $j=1, \ldots ,
n-\ell-k_i+1$) be a collection of complex analytic 1-forms on a neighbourhood of the origin
in
$\CC^n$ the restriction of which has no singular points on $V$ outside of the origin in a
neighbourhood of it. Let $I_{V,\{ \omega_j^{(i)}\}}$ be the
ideal in the ring
${\cal O}_{\CC^n,0}$ generated by the functions $f_1, \ldots , f_\ell$ and by the
$(n-k_i+1) \times (n-k_i+1)$ minors of all the matrices
$$(df_1(x), \ldots , df_\ell(x), \omega_1^{(i)}(x), \ldots ,
\omega_{n-\ell-k_i+1}^{(i)}(x))$$
for all $i=1, \ldots, s$.

\begin{theorem} \label{Th4}
$${\rm ind}_{V,0}\{ \omega^{(i)}_j \} = \dim_\CC {\cal O}_{\CC^n,0}/I_{V,\{
\omega_j^{(i)}\}}.$$
\end{theorem}

See the ''proof'' of Theorem~\ref{Th2}.

\addvspace{3mm}

\noindent {\bf Remark.} Let $\{ X_j^{(i)} \}$ be a collection of vector fields on a
neighbourhood of the origin in $(\CC^n,0)$ ($i=1, \ldots , s$; $j=1, \ldots , n-\ell-k_i+1$;
$\sum k_i = n-\ell$) which are tangent to the ICIS $(V,0)=\{f_1= \cdots = f_\ell=0\}
\subset (\CC^n,0)$ at non-singular points of $V$. One can define the index ${\rm ind}_{V,0}
\{ X_j^{(i)} \}$ as the degree of the mapping $K \to W_{n, \widehat{\bf k}}$ which sends a
point $x \in K$ to the collection of $n \times (n-k_i+1)$ matrices
$$\{ ( \mbox{grad}\, f_1(x), \ldots , \mbox{grad}\, f_\ell(x), X_1^{(i)}(x), \ldots ,
X_{n-\ell-k_i+1}^{(i)}(x))\}, \quad i=1, \ldots , s.$$
Here 
$$\mbox{grad}\, f_r(x)=\left( \overline{\frac{\partial f_r}{\partial x_1}(x)}, \ldots, 
\overline{\frac{\partial f_r}{\partial x_n}(x)} \right),$$
where $\overline{z}$ is the complex conjugate of the complex number $z$. 
For $s=1$, $k_1=n-\ell$, this definition coincides with the definition of the index of a
vector field on an ICIS from \cite{GSV, SS}.
For vector fields the
analogue of Theorem~\ref{Th3} holds with the only difference that the sum of the indices is
equal to the characteristic number
$\langle \prod_{i=1}^s
c_{k_i}(T\widetilde{V} \otimes L), [ \widetilde{V} ] \rangle$. However, a formula similar
to that of Theorem~\ref{Th4} does not exist.  A reason for that is that in this case 
the index is the intersection number
with $D_{n, \widehat{\bf k}}$ of the image of the ICIS $(V,0)$ under a map which is not
complex analytic. Moreover, in some cases this index can be
negative.

\bigskip
\noindent Universit\"{a}t Hannover, Institut f\"{u}r Mathematik \\
Postfach 6009, D-30060 Hannover, Germany \\
E-mail: ebeling@math.uni-hannover.de\\

\medskip
\noindent Moscow State University, Faculty of Mechanics and Mathematics\\
Moscow, 119992, Russia\\
E-mail: sabir@mccme.ru

\end{document}